\documentclass{ifacconf}

\usepackage{graphicx}      
\usepackage[square,numbers]{natbib}        
\usepackage{tabularx}
\usepackage{multirow}
\usepackage{rotating}
\usepackage{csquotes}
\usepackage{amsfonts}
\usepackage{xfrac}
\usepackage{empheq}
\usepackage{fancyhdr}
\usepackage{mathtools}
\usepackage{inputenc}
\usepackage{url}
\usepackage{fix-cm}
\usepackage{bm}
\usepackage{color}
\usepackage{xcolor}
\usepackage{etoolbox}
\usepackage{latexsym}
\usepackage{textcomp}
\usepackage{setspace}
\DeclareGraphicsExtensions{.eps,.pdf,.jpg}
\usepackage{epstopdf}
\theoremstyle{plain}
\usepackage{lipsum}

\usepackage{bibentry}

\newtheorem{assumption}{First-Order Approximation}
\newtheorem{proposition}{Proposition}
\newtheorem{lemma}{Lemma}

\newtheorem{theo}{Theorem}

\newcommand{\vect}[1]{\mathbf{#1}}
\newcommand{\bxi}{\boldsymbol{\xi}}

\begin{document}
\begin{frontmatter}

\title{Analytical Initialization of a Continuation-Based Indirect Method for Optimal Control of Endo-Atmospheric Launch Vehicle Systems} 

\author[First,Second]{Riccardo Bonalli} 
\author[Second]{Bruno H\'{e}riss\'{e}} 
\author[First]{Emmanuel Tr\'{e}lat}

\address[First]{Sorbonne Universit\'{e}s, UPMC Univ Paris 06, CNRS UMR 7598, Laboratoire Jacques-Louis Lions, F-75005, Paris, France.}
\address[Second]{Onera - The French Aerospace Lab, F~-~91761 Palaiseau, France. \begingroup \footnotesize riccardo.bonalli@onera.fr , bruno.herisse@onera.fr , emmanuel.trelat@upmc.fr \endgroup}

\begin{abstract}                
In this paper, we propose a strategy to solve endo-atmospheric launch vehicle optimal control problems using indirect methods. More specifically, we combine shooting methods with an adequate continuation algorithm, taking advantage of the knowledge of an analytical solution of a simpler problem. This procedure is resumed in two main steps. We first simplify the physical dynamics to obtain an analytical guidance law which is used as initial guess for a shooting method. Then, a continuation procedure makes the problem converge to the complete dynamics leading to the optimal solution of the original system. Numerical results are presented.
\end{abstract}


\end{frontmatter}

\section{Introduction}

Guidance of autonomous launch vehicles towards rendezvous points is a complex task often considered in missile applications, mainly for interception of targets. It represents an optimal control problem, whose aim consists in finding a control law enabling the vehicle to join a final point of the 3D space considering prescribed constraints as well as performance criteria. The rendezvous point may be a static point, for example an intermediate point at the middle of the flight mission, as well as a moving point if, for example, the task consists in intercepting a maneuvering target. This requires not only \textit{high numerical precision} of concerned algorithms but also a \textit{real-time processing} of optimal trajectories.

The most common approach to solve this kind of task resides on \textit{analytical guidance laws}. They correct errors coming from perturbations and misreading of the system. For interceptor missile applications, in \citep{LinTsai}, \citep{Lin}, \citep{Nahshon} guidance laws that maximize the final velocity and so the lethality of the interceptor are designed. Varying velocity is also considered for the design of an energy optimal guidance law in \citep{Taek}. However, the trajectories induced by guidance laws are not optimal because of some approximations made.

Ensuring the optimality of trajectories can be achieved exploiting \textit{direct methods}. These techniques consist in discretizing each component of the optimal control problem (the state, the control, etc.) to reduce the whole mathematical representation to a nonlinear constrained optimization problem. Since they are quite robust, they are widely used \citep{Paris}, \citep{RossOld}, \citep{Ross}. However, these methods are computationally demanding and can often be used offline uniquely.

In order to manage efficiently real-time processing of optimal control sequences for launch vehicle systems one may consider \textit{indirect methods}. These use a mathematical study of the system (exploiting the \textit{Pontryagin Maximum Principle}) to determine some necessary conditions of optimality. Indirect methods converge much faster than direct methods with a better precision. Since the problem is equivalent to the research of zeros of a function, the main difficulty remains their \textit{initialization}. In \citep{Lu1}, \citep{Lu3} the initialization problem is bypassed using finite differences algorithms and multiple-shooting methods respectively. By contrast, in \citep{Pontani} the authors use second-order conditions and conjugate point theory for a multistage launch vehicle problem; the initialization of the numerical method is based on the particle swarm algorithm. However, most of these approaches remain computationally demanding and not easily applicable in view of real-time processing.

In this paper, we propose to solve endo-atmospheric launch vehicle optimal control problems using indirect methods managing the issue coming from the initialization by both an analytical guidance law and a \textit{continuation method}. Continuation procedures have shown to be reliable and robust for problems like coplanar orbit transfer and atmospheric reentry \citep{EmmanuelH}, \citep{Petit}. This combination allows to preserve precision and fast numerical computations. The proposed approach is resumed in two main steps. We first simplify the physical dynamics to obtain an analytical guidance law which is used as initial guess for a shooting method. Then, a continuation procedure makes the problem converge to the complete dynamics leading to the optimal solution of the original system.

The paper is structured as follows. Section 2 introduces the dynamics of a general endo-atmospheric launch vehicle system, the optimal control formulation and the related numerical approach. Section 3 is devoted to the construction of a simplified dynamics able to initialize successfully a shooting method and to make the final algorithm converge efficiently to the complete dynamics. In Section 4 the continuation method to recover the original dynamics is presented with its indirect formulation and numerical tests. Finally, Section 5 contains conclusions and perspectives.


\section{Dynamics, Optimal Control Problem and Numerical Method}

\subsection{Physical Model}

Let $(\mathbf{I}, \mathbf{J}, \mathbf{K})$ be an inertial frame centered at the center of the planet $O$, $(\mathbf{e}_L, \mathbf{e}_l, \mathbf{e}_r)$ be the NED frame and $(\mathbf{i}, \mathbf{j}, \mathbf{k})$ be the velocity frame. The endo-atmospheric launch vehicle system is modeled as an axisymmetric thrust propelled rigid body of mass $m$ \citep{brunoMaster}. The coordinates $(r, l, L, v, \gamma, \chi)~\in~\mathbb{R}^6$ ($L$ is the latitude, $l$ is the longitude, $\gamma$ is the path angle and $\chi$ is the azimuth angle) are used to represent the position $\bxi = (r \cos(L) \cos(l), r \cos(L) \sin(l), r \sin(L))$ of the center of mass $G$ of the vehicle and its velocity $\mathbf{v}~=~v \cos(\gamma) \cos(\chi) \mathbf{e}_L + v \cos(\gamma) \sin(\chi) \mathbf{e}_l - v \sin(\gamma) \mathbf{e}_r$.

Neglecting the wind velocity, the Coriolis and the centripetal force (this is legitimate because of the short length of the considered trajectories), the dynamics takes the following form
\begingroup
\footnotesize
\begin{eqnarray} \label{dynamics}
\dot{r} &=& v\sin(\gamma) \; , \qquad \dot{L} = \displaystyle \frac{v}{r} \cos(\gamma) \cos(\chi) \; , \qquad \dot{l} = \displaystyle \frac{v}{r} \frac{\cos(\gamma) \sin(\chi)}{\cos(L)} \medskip \nonumber \\
\dot{v} &=& \displaystyle \frac{f_T}{m} \cos(\alpha) - (d + \eta c_m u^2) v^2 - g \sin(\gamma) \medskip \nonumber \\
\dot{\gamma} &=& \displaystyle \frac{f_T}{m v} \sin(\alpha) \cos(\beta) + v c_m u \cos(\beta) + \left(\frac{v}{r} - \frac{g}{v}\right) \cos(\gamma) \medskip \\
\dot{\chi} &=& \displaystyle \frac{f_T}{m v} \frac{\sin(\alpha)}{\cos(\gamma)} \sin(\beta) + \frac{v c_m}{\cos(\gamma)} u \sin(\beta) + \frac{v}{r} \cos(\gamma) \sin(\chi) \tan(L) \medskip \nonumber \\
\dot{m} &=& -q \nonumber
\end{eqnarray}
\endgroup
where $g$ is the modulus of the gravity, $\eta$ is an aerodynamic efficiency factor, $\alpha = \alpha_{\max} u$ is the \textit{angle of attack} while $\beta$ is the \textit{angle of bank}, $u$ stands for the normalized lift coefficient while $q = q(t)$ is the mass flow and $f_T = f_T(t)$ represents the modulus of the propulsion which depends directly on $q(t)$.

Based on a standard model of flight dynamics \citep{brunoMaster}, \citep{Bruno}, coefficients $d$ and $c_m$ are approximated by $d~=~d(r)~~=~\frac{1}{2 m} \rho SC_{D_0}$ and $c_m = c_m(r) = \frac{1}{2 m} \rho SC_{L_{\max}}$ where $S$ is the reference area, $C_{L_{\max}}$ is the maximal value of the lift coefficient and $C_{D_0}$ is the drag coefficient for $\alpha = 0$, which is considered constant; finally, $\rho$ is the air density for which an exponential model $\rho~=~\rho(r)~=~\rho_0 \exp(-(r~-~r_T)/h_r)$ is considered, where $h_r$ is a reference altitude. Since $q(t)$ is a predefined function of time, in this paper controls are only $u$ and $\beta$.


\subsection{Optimal Control Problem and Maximum Principle}

%
Consider now the Optimal Control Problem (\textbf{OCP})
\begin{eqnarray} \label{generalDyn}
\displaystyle \min \  && \hspace{-0.2cm} \int_0^{t_f} f^0(t,\vect{x}(t),\vect{u}(t)) \; dt \medskip \nonumber \\
\dot{\vect{x}}(t) &=& \vect{f}(t,\vect{x}(t),\vect{u}(t)) \medskip \\
\vect{x}(0) &=& \vect{x}_0 \in M_0 \; , \ \ \vect{x}(t_f) = \vect{x}_{t_f} \in M_f \nonumber
\end{eqnarray}
where $\vect{x}(t) = (r,L,l,v,\gamma,\chi)(t) \in \mathbb{R}^6$, $\vect{u}(t) = (u,\beta)(t) \in \mathbb{R}^2$, $\vect{f}$ is the mapping defined by dynamics (\ref{dynamics}), $M_0$, $M_f$ are smooth submanifolds in $\mathbb{R}^6$ and the transfer time $t_f$ is not fixed. Finally,
\begin{equation} \label{cost}
f^0 = \displaystyle \sigma u^2 - \bigg( \frac{f_T}{m} \cos(\alpha) - (d + \eta c_m u^2) v^2 - g \sin(\gamma) \bigg)
\end{equation}
where $\sigma \geq 0$ is constant. By definition, $u^2$ takes its values in $[0,1]$. However, we do not consider any boundaries on $u$ preferring to penalize it, using $\sigma$, $\eta c_m v^2$ within the cost.

The \textit{Pontryagin Maximum Principle} (PMP) \citep{pontryagin} states that, if $\vect{u}$ is optimal with response defined on $[0,t_f]$, and shortly denoted $\vect{x}(t)$, then there exists $p^0 \leq 0$ and $\vect{p}~\in~AC([0,t_f],\mathbb{R}^6)$ such that a.e. on $[0,t_f]$
\begingroup
\begin{eqnarray} \label{primalDual}
\displaystyle &\dot{\vect{x}} = \frac{\partial H}{\partial \vect{p}} \; , \ \dot{\vect{p}} = -\displaystyle \frac{\partial H}{\partial \vect{x}} \; , \ H(t,\vect{x},\vect{p},\vect{u}) = \max_{\vect{v} \in U} H(t,\vect{x},\vect{p},\vect{v}) \nonumber \medskip \\
&\max_{\vect{v} \in U} H(t_f,\vect{x}(t_f),\vect{p}(t_f),\vect{v}) = 0
\end{eqnarray}
\endgroup
where $H = \left \langle \vect{p} , \vect{f}(t,\vect{x},\vect{u}) \right \rangle + p^0 f^0(t,\vect{x},\vect{u})$ is the Hamiltonian and $\vect{p}$ satisfies the transversality conditions
\begin{equation} \label{transv}
\vect{p}(0) \perp T_{\vect{x}_0} M_0 \; , \quad \vect{p}(t_f) \perp T_{\vect{x}_{t_f}} M_f 
\end{equation}

Treating (\textbf{OCP}) by indirect methods consists in solving
\begin{eqnarray} \label{adjointSys}
\dot{\vect{x}}(t) &=& \vect{f}(t,\vect{x}(t),\vect{p}(t)) \; , \quad \ \, \vect{x}(0) = \vect{x}_0 \; , \ \vect{x}(t_f) = \vect{x}_f \\
\dot{\vect{p}}(t) &=& -\frac{\partial H}{\partial \vect{x}}(t,\vect{x}(t),\vect{p}(t)) \; , \quad \vect{p}(0) = \vect{p}_0 \nonumber
\end{eqnarray}
with an appropriate value for $\vect{p}_0$.

\subsection{Shooting and Continuation Method}

It is known \citep{Emmanuel} that finding a solution of (\ref{adjointSys}) can be reduced to solve $G(\vect{p}_0,t_f) = 0 \; , \ G : \mathbb{R}^{n+1} \longrightarrow \mathbb{R}^{n+1}$ ($G$ is called shooting function) using Newton-type methods. This is the content of the well known \textit{shooting method} in optimal control \citep{Emmanuel}, \citep{EmmanuelReview}, \citep{Stoer}. Its advantage is its extremely good numerical accuracy, relevant for aerospace applications \citep{EmmanuelReview}. Since it relies on the Newton method, it inherits of the very quick convergence properties of the Newton method. The main drawback of the shooting method is that it may be difficult to initialize.

To overcome this difficulty, one can entrust with the robustness of the \textit{continuation method}. It consists in deforming the problem into a simpler one that we are able to solve and then in solving a series of shooting problems, step by step by parameter deformation, to recover the original problem \citep{Werner}. This approach increases the efficiency of the shooting method because it allows to relax its initialization. The continuation parameter $\lambda$ may be a physical parameter (or several) of the problem, or an artificial one. The path consists of a convex combination of the simpler problem and of the original one, with $\lambda \in [0,1]$.

The main algorithm consists then in finding a solution of some simplification of (\textbf{OCP}) first and, from this, solving by continuation the original formulation (\textbf{OCP}).

\section{New Guidance Law as Good Estimate for Continuation Method}

Continuation methods allow us to solve iteratively (\textbf{OCP}) once the solution of some (usually) simpler optimal control problem is known. Here, we introduce and treat an efficient simpler problem coming from a modification of (\ref{dynamics}).

This modified version of (\textbf{OCP}) is designed with the hope that, on one hand, the shooting method can be easily initialized and, on the other hand, the continuation is feasible. It is at this stage that intuition of aerospace engineer is of primary importance, in order to design simplified problems and homotopy parameters resulting in a meaningful continuation procedure.

\subsection{Simplified Problem and Analytical Smooth Controls}

If one ignores the contributions of the curvature of the Earth, of the gravity and of the propulsion within (\ref{dynamics}), introducing the curvilinear abscissa $s(t) := \int_{0}^{t} v(t') \; dt'$ and a new variable $w := \ln (v)$, the following simplification of (\textbf{OCP}) is obtained (the equation of $w'$ can be neglected because it does not influence anymore the dynamics)
\begingroup
\footnotesize
\begin{eqnarray} \label{dynamicsSimplyAbs}
\quad \displaystyle &\min& \int_{0}^{s_f} (d + \eta c_m (u_1^2 + u_2^2)) \; dt \medskip \nonumber \\
&r'& = \displaystyle \frac{d r}{d s} = \sin(\gamma) \; , \qquad L' = \displaystyle \frac{d L}{d s} = \displaystyle \frac{\cos(\gamma) \cos(\chi)}{r} \medskip \nonumber \\
&l'& = \displaystyle \frac{d l}{d s} = \displaystyle \frac{\cos(\gamma) \sin(\chi)}{r \cos(L)} \medskip \\
&\gamma'& = \displaystyle \frac{d \gamma}{d s} = c_m u_1 \; , \qquad \chi' = \displaystyle \frac{d \chi}{d s} = \displaystyle \frac{c_m u_2}{\cos(\gamma)} \nonumber
\end{eqnarray}
\endgroup
where $u_1 := u \cos (\beta)$, $u_2 := u \sin (\beta)$ and with initial, final conditions $(r_0,L_0,l_0,v_0,\gamma_0,\chi_0)$, $(r_f,L_f,l_f,\gamma_f,\chi_f)$.

The next step consists in finding an analytical solution of (\ref{dynamicsSimplyAbs}). This is achieved exploiting the PMP formulation, as showed in the following.

The associated Hamiltonian function is
\begingroup
\footnotesize
\begin{equation} \label{hamSimply} \begin{array}{c} 
\displaystyle H = p_r \sin (\gamma) + p_L \frac{\cos(\gamma) \cos(\chi)}{r} + p_l \frac{\cos(\gamma) \sin(\chi)}{r \cos(L)} \medskip \\
\displaystyle + p_{\gamma} c_m u_1 + p_{\chi} \frac{c_m u_2}{\cos(\gamma)} + p^0 (d + \eta c_m (u_1^2 + u_2^2))
\end{array} 
\end{equation}
\endgroup

It can be remarked that, since the transfer time is not fixed and formulation (\ref{dynamicsSimplyAbs}) is autonomous, (\ref{hamSimply}) takes zero as value for all times $t \in [0,t_f]$ \citep{pontryagin}.
Since no constraints are considered on $u_1$, $u_2$, from applying the weaker version of PMP \citep{pontryagin}, it follows explicitly that
\begin{equation} \label{controlExp}
p_{\gamma} = -2 \eta p^0 u_1 \; , \quad p_{\chi} = -2 \eta p^0 \cos(\gamma) u_2
\end{equation}

Simple calculations show that
\begingroup
\scriptsize
\begin{eqnarray} \label{dualSys}
\displaystyle &p_r'& = p_L \frac{\cos(\gamma) \cos(\chi)}{r^2} + p_l \frac{\cos(\gamma) \sin(\gamma)}{r^2 \cos(L)} + \frac{p^0}{h_r} (d - \eta c_m (u_1^2 + u_2^2)) \medskip \nonumber \\
\displaystyle &p_L'& = -p_l \frac{\cos(\gamma) \sin(\chi) \tan(L)}{r \cos(L)} \; , \qquad p_l' = 0 \medskip \\
\displaystyle &p_{\gamma}'& = -p_r \cos(\gamma) + p_L \frac{\sin(\gamma) \cos(\chi)}{r} + p_l \frac{\sin(\gamma) \sin(\chi)}{r \cos(L)} - p_{\chi} \frac{c_m u_2 \tan(\gamma)}{\cos(\gamma)} \medskip \nonumber \\
\displaystyle &p_{\chi}'& = p_L \frac{\cos(\gamma) \sin(\chi)}{r} - p_l \frac{\cos(\gamma) \cos(\chi)}{r \cos(L)} \nonumber
\end{eqnarray}
\endgroup

A further analysis streamlines this formulation letting one considering only \textit{normal trajectories}.

\begin{proposition} \label{prop1}
\textit{Assume that $\cos(\gamma) \neq 0$ a.e. in $[0,t_f]$. Then, no abnormal trajectories arise in (\ref{dynamicsSimplyAbs}), i.e. $p^0 \neq 0$.}
\end{proposition}

\begin{pf} If $p^0 = 0$, then from (\ref{controlExp}) it follows $p_{\gamma} = p_{\chi} \equiv 0$ in $[0,t_f]$. Using $H = 0$, $p_{\gamma}' = 0$ and $p_{\chi}' = 0$ raises a homogeneous linear system of equations in $(p_r,p_L,p_l)$ whose associated matrix has $-\frac{\cos(\gamma)}{r^2 \cos(L)}$ as determinant which is not zero almost everywhere in $[0,t_f]$ by assumption. Then $(p_r(t),p_L(t),p_l(t),p_{\gamma}(t),p_{\chi}(t),p^0) \equiv 0$ on $[0,t_f]$ which raises a contradiction with the PMP. $_\Box$ \\
\end{pf}

The assumption $\cos(\gamma) \neq 0$ a.e. in $[0,t_f]$ is relevant because, if the trajectory takes $\cos(\gamma) = 0$ as optimal value, a change of coordinates like $\tilde{\gamma} = \gamma + \frac{\pi}{2}$ can be done on the dynamics during the numerical integration to avoid this singularity. Moreover, other changes of coordinates that optimize the computational speed can be taken into account \citep{NewEmmanuel}. Thus, by Proposition \ref{prop1} and the PMP, we can substitute $p^0 = -1$ within the previous simplified formulation.

Proceeding formally under smoothness assumption for optimal controls of this simplified problem, analyzing the second derivative of $u_1$ and $u_2$ allows to seek some relations. Indeed, using (\ref{controlExp}) and (\ref{dualSys}) it can be easily shown that
\begin{footnotesize}
\begin{eqnarray*}
\displaystyle u_1'' = \frac{p_{\gamma}''}{2 \eta} = \frac{d}{2 \eta} \left( c_m u_1 + \frac{\cos(\gamma)}{h_r} \right) - \frac{1}{2 \eta r^2} \left( p_L \cos(\chi) + p_l \frac{\sin(\chi)}{\cos(L)} \right) \medskip \\
+ P_1 (u_1,u_2,u_1',u_2')
\end{eqnarray*}
\end{footnotesize}
\hspace{-0.42cm} where $P_1 (u_1,u_2,u_1',u_2')$ is a polynomial of degree greater than one of $u_1$, $u_2$ and their first derivatives. Since $\frac{1}{r^2} \ll 1$ the previous expression can be approximated. Iterating the same procedure on $u_2$ one is led to solve ($P_2 (u_1,u_2,u_1',u_2')$ identifies again a polynomial of degree greater than one of $u_1$, $u_2$ and their first derivatives)
\begin{eqnarray} \label{u1u2App}
\begin{cases}
\displaystyle u_1'' = \frac{c_m d}{2 \eta} \left( u_1 + \frac{\cos(\gamma)}{c_m h_r} \right) + P_1 (u_1,u_2,u_1',u_2') \medskip \\
\displaystyle u_2'' = \frac{c_m d}{2 \eta} u_2 + P_2 (u_1,u_2,u_1',u_2')
\end{cases}
\end{eqnarray}

It is clear that, if additional assumptions allow to neglect the contribution of $P_1 (u_1,u_2,u_1',u_2')$, $P_2 (u_1,u_2,u_1',u_2')$, equations of system (\ref{u1u2App}) can be solved independently using the assumption that parameters $c_m$, $d$ and $\cos(\gamma)$ are almost constant along the trajectory, justifying this approximation because scenarios with limited altitude are considered to solve the simplified problem (\ref{dynamicsSimplyAbs}). This is not limiting, because, once the shooting method on simplified problem has converged, a continuation procedure on the final point can be started changing the initial scenario.

\subsection{Local Controllability of the Simplified System}

An evident sufficient condition that allows to neglect the contribution of $P_1 (u_1,u_2,u_1',u_2')$, $P_2 (u_1,u_2,u_1',u_2')$ within system (\ref{u1u2App}) is to deal with small (in the sense of $L^{\infty}$ norm) controls $\vect{z} := (u_1,u_2)$ that have small derivatives. This is achieved studying the \textit{local controllability around control} $\bar{\vect{z}} = (0,0)$ in
$$
W^{1,\infty}(0,t_f;\mathbb{R}^p) := \{ \vect{z} \in L^{\infty}(0,t_f;\mathbb{R}^p) : \vect{z}' \in L^{\infty}(0,t_f;\mathbb{R}^p) \}
$$

where $\vect{z}'$ is the first distributional derivative of $\vect{z}$. The set of admissible controls on $[0, t_f]$ is denoted by $\mathcal{U}_{\vect{y}_0, t_f, \mathbb{R}^p}$ where $\vect{y}_0$~denotes the initial point of the state variable $\vect{y}~:=~(r,L,l,\gamma,\chi)$.

The attention is focused on a special class of control system to which the dynamics of (\ref{dynamicsSimplyAbs}) belongs, the \textit{affine systems}. A control system (\ref{generalDyn}) is affine if
\begin{equation} \label{affine}
\vect{f} = \vect{f} (\vect{y},\vect{z}) = \vect{f}_0(\vect{y}) + \sum_{i=1}^m z_i \vect{f}_i(\vect{y})
\end{equation}
where the $\vect{f}_i$'s are (at least) $C^1$ vector fields. This standard result holds (whose proof will be reported in \citep{mio})

\begin{theo} \label{contLoc}
Let $\vect{y}_0 \in \mathbb{R}^n$ and $\bar{\vect{z}} \in \mathcal{U}_{\vect{y}_0, t_f, \mathbb{R}^p}$. If the linearized system of (\ref{affine}) along $(\vect{y}_{\bar{\vect{z}}}(\cdot),\bar{\vect{z}}(\cdot))$ is controllable from $\vect{y}_0$ in time $t_f$, then (\ref{affine}) is locally controllable at $\vect{y}_{\bar{\vect{z}}}(t_f)$ from $\vect{y}_0$ in time $t_f$ using controls $\vect{z}$ with $\| \vect{z} - \bar{\vect{z}} \|_{W^{1,\infty}}$ small enough.
\end{theo} 

All these results can be applied to the dynamics of (\ref{dynamicsSimplyAbs}) allowing to give conditions to solve easily (\ref{u1u2App}). The idea is to compute the trajectory solution of (\ref{dynamicsSimplyAbs}) given by control $\bar{\vect{z}}~=~(0,0) \in~W^{1,\infty}(0,t_f;\mathbb{R}^2)$ and to consider scenarios with final point close to $\vect{y}_{\bar{\vect{z}}}(t_f)$. \textit{If the linearized system along} $(\vect{y}_{\bar{\vect{z}}}(\cdot),\bar{\vect{z}}(\cdot))$ \textit{is controllable from} $\vect{y}_0$ \textit{in time} $t_f$, then, thanks to Theorem~\ref{contLoc}, the dynamics of (\ref{dynamicsSimplyAbs}) is locally controllable at $\vect{y}_{\bar{\vect{z}}}(t_f)$ from $\vect{y}_0$ in time $t_f$ with small controls in $W^{1,\infty}(0,t_f;\mathbb{R}^2)$. In other words, if the considered scenarios have final point close to $\vect{y}_{\bar{\vect{z}}}(t_f)$ we can neglect the contribution of $P_1 (u_1,u_2,u_1',u_2')$, $P_2 (u_1,u_2,u_1',u_2')$ within system (\ref{u1u2App}) to solve it analytically. It must be noted that it is not restrictive considering scenarios with final point close to $\vect{y}_{\bar{\vect{z}}}(t_f)$ to solve (\ref{dynamicsSimplyAbs}), because, as said before, a continuation procedure on the final point can be started changing considerably the initial scenario.

The linearized system of (\ref{dynamicsSimplyAbs}) has the form
$$
\vect{y}'(s) = A(s) \vect{y}(s) + B(s) \vect{z}(s)
$$
where $A(s)$ is the matrix
$$
\tiny \hspace{-0.2cm} \left( \begin{array}{ccccc}
0 & 0 & 0 & \cos(\gamma_0) & 0 \\
-\frac{\cos(\gamma_0) \cos(\chi_0)}{r^2} & 0 & 0 & -\frac{\sin(\gamma_0) \cos(\chi_0)}{r} & -\frac{\cos(\gamma_0) \sin(\chi_0)}{r} \\
-\frac{\cos(\gamma_0) \sin(\chi_0)}{r^2 \cos(L)} & \frac{\cos(\gamma_0) \sin(\chi_0) \tan(L)}{r \cos(L)} & 0 & -\frac{\sin(\gamma_0) \sin(\chi_0)}{r \cos(L)} & \frac{\cos(\gamma_0) \cos(\chi_0)}{r \cos(L)} \\
0 & 0 & 0 & 0 & 0 \\
0 & 0 & 0 & 0 & 0
\end{array} \right)
$$
and
$$
\scriptsize
B(s) = \footnotesize \left( \begin{array}{ccccc}
0 & 0 & 0 & c_m(r) & 0\\
0 & 0 & 0 & 0 & \frac{c_m(r)}{\cos(\gamma_0)}\\
\end{array} \right)^{\top}
$$

This system is controllable in any $s_f$. To show this, we consider the well known fact \citep{Emmanuel} that, for a system with no constraint $\vect{y}'(s) = A(s) \vect{y}(s) + B(s) \vect{z}(s)$ where $s~\mapsto~A(s)$ and $s \mapsto B(s)$ are of class $C^{\infty}$, if, defined the sequence $B_0(s)~:=~B(s)$, $B_{i+1}(s) := A(s) B_i(s) - \frac{d B_i}{ds}(s)$ for $i \in \mathbb{N}$, there exists $s \in [0,s_f]$ such that
$$
\textnormal{Span} \{ B_i(s) \vect{w} : \vect{w} \in \mathbb{R}^p \; , \; i \in \mathbb{N} \} = \mathbb{R}^n
$$
then the system is controllable in time $s_f$. Denoting $\vect{w}_1~=~(1,0)^{\top}$, $\vect{w}_2 = (0,1)^{\top}$, it is easy to show that
\begingroup
\footnotesize
$$
\textnormal{det} \big( B_0(s) \vect{w}_1 , B_0(s) \vect{w}_2 , B_1(s) \vect{w}_1 , B_1(s) \vect{w}_2 , B_2(s) \vect{w}_1 \big) = \frac{c_m^5}{r^3 \cos(L)}
$$
\endgroup
which is different from zero. Then, the claim follows.

This result validates the previous approach. Initial conditions for controls and their first derivatives must be sought to solve (\ref{u1u2App}) (where we impose then $P_1 = P_2 = 0$). Two conditions are given using the two last equations of (\ref{dynamicsSimplyAbs}). The second ones can be investigated exploiting first an analysis of the \textit{Line Of Sight} (LOS).

\subsection{LOS Analysis}

\vspace{-0.1cm}

The line of sight is the vector joining the current position $\bxi$ to the desired final point $\bxi_f$. It results to be useful defining its modulus $R~:=\| \bxi_f - \bxi \|$ and the associated unitary vector $\vect{n}~:=~\frac{\bxi_f - \bxi}{R}$. We denote $(\mathbf{i}_{n}, \mathbf{j}_{n}, \mathbf{k}_{n})$ the orthonormal frame where $\mathbf{i}_{n} = \vect{n}$, $\mathbf{j}_{n}$ is the unitary vector in the plane $(\mathbf{i}_{n}, \mathbf{e}_r)$ perpendicular to $\mathbf{i}_{n}$ and oriented by $\mathbf{j}_{n} \cdot \mathbf{e}_r < 0$ and $\mathbf{k}_{n} = \mathbf{i}_{n} \wedge~\mathbf{j}_{n}$. Frames $(\mathbf{i},\mathbf{j},\mathbf{k})$, $(\mathbf{i}_{n},\mathbf{j}_{n},\mathbf{k}_{n})$ are functions of $(\mathbf{e}_L, \mathbf{e}_l, \mathbf{e}_r)$ by definition \citep{CNES}, locally according to

\begingroup
\begin{eqnarray*}
&\mathbf{i}& = \frac{\vect{v}}{v} = \cos(\gamma) \cos(\chi) \mathbf{e}_L + \cos(\gamma) \sin(\chi) \mathbf{e}_l - \sin(\gamma) \mathbf{e}_r \\
&\mathbf{j}& = -\sin(\gamma) \cos(\chi) \mathbf{e}_L - \sin(\gamma) \sin(\chi) \mathbf{e}_l - \cos(\gamma) \mathbf{e}_r \\
&\mathbf{k}& = -\sin(\chi) \mathbf{e}_L + \cos(\chi) \mathbf{e}_l \\
\hspace{-0.1cm} &\mathbf{i}_{n}& \hspace{-0.05cm} = \vect{n} = \cos(\lambda_1) \cos(\lambda_2) \mathbf{e}_L + \cos(\lambda_1) \sin(\lambda_2) \mathbf{e}_l - \sin(\lambda_1) \mathbf{e}_r \\
&\mathbf{j}_{n}& = -\sin(\lambda_1) \cos(\lambda_2) \mathbf{e}_L - \sin(\lambda_1) \sin(\lambda_2) \mathbf{e}_l - \cos(\lambda_1) \mathbf{e}_r \\
&\mathbf{k}_{n}& = -\sin(\lambda_2) \mathbf{e}_L + \cos(\lambda_2) \mathbf{e}_l
\end{eqnarray*}
\endgroup

where $\lambda_1$, $\lambda_2$ are Euler angles. We denote $M_{(\mathbf{i},\mathbf{j},\mathbf{k})}$, $M_{(\mathbf{i}_{n},\mathbf{j}_{n},\mathbf{k}_{n})}$ the associated transition matrices from frame $(\mathbf{e}_L, \mathbf{e}_l, \mathbf{e}_r)$ to frame $(\mathbf{i},\mathbf{j},\mathbf{k})$ and from frame $(\mathbf{e}_L, \mathbf{e}_l, \mathbf{e}_r)$ to frame $(\mathbf{i}_{n},\mathbf{j}_{n},\mathbf{k}_{n})$ respectively. The coordinates of a vector within a frame $(\mathbf{a},\mathbf{b},\mathbf{c})$ are denoted as $( \cdot )_{(\mathbf{a},\mathbf{b},\mathbf{c})}$. The following analysis seeks some approximate relations between the derivatives of $\lambda_1$, $\lambda_2$ and angles $\gamma$, $\chi$ under appropriate assumptions.

\vspace{0.4cm}

\begin{assumption}[\textbf{FOA}]
Along the considered trajectories, only weak variations from the LOS direction are allowed, i.e. $\gamma \cong \lambda_1$~, $\chi \cong \lambda_2$. Moreover, the contribution of the rotation of frame $(\mathbf{e}_L, \mathbf{e}_l, \mathbf{e}_r)$ around frame $(\mathbf{I},\mathbf{J},\mathbf{K})$ is neglected, i.e. $\frac{d}{dt}( \vect{n} )_{(\mathbf{e}_L, \mathbf{e}_l, \mathbf{e}_r)} \cong \left( \frac{d}{dt} \vect{n} \right)_{(\mathbf{e}_L, \mathbf{e}_l, \mathbf{e}_r)}$.
\end{assumption}

It must be remarked that this assumption is not a strong one because the initial scenario involved lives in a neighborhood of the null control and the considered trajectories are short enough to make this choice effective. It allows to manipulate transition maps obtaining easily the desired results.

\newpage

\begin{proposition} \label{PropNeeded}
Under (\textnormal{\textbf{FOA}}), the transition matrix $M$ from frame $(\mathbf{i},\mathbf{j},\mathbf{k})$ to frame $(\mathbf{i}_{n},\mathbf{j}_{n},\mathbf{k}_{n})$ is
\[\footnotesize \left( \begin{array}{ccc}
\cos(\gamma - \lambda_1) & -\sin(\gamma - \lambda_1) & -\sin(\chi - \lambda_2) \cos(\gamma) \\
\sin(\gamma - \lambda_1) & \cos(\gamma-\lambda_1) & \sin(\chi - \lambda_2) \sin(\gamma) \\
\sin(\chi - \lambda_2) \cos(\gamma) & -\sin(\chi - \lambda_2) \sin(\gamma) & \cos(\chi - \lambda_2) \end{array} \right) \]
Moreover, we have
\[(\mathbf{v})_{(\mathbf{i}_{n},\mathbf{j}_{n},\mathbf{k}_{n})} = \left(\footnotesize \begin{array}{c} v \cos(\gamma - \lambda_1) \\ v \sin(\gamma - \lambda_1) \\
v \sin(\chi - \lambda_2) \cos(\gamma) \\ \end{array} \right) \]
\end{proposition}

\begin{pf} Matrix $M$ is obtained easily by approximating $\cos(\chi - \lambda_2) \cong 1$ within the changing variable rotation $M_{(\mathbf{i}_{n},\mathbf{j}_{n},\mathbf{k}_{n})} \circ M_{(\mathbf{i},\mathbf{j},\mathbf{k})}^{\top}$. The second relation follows by multiplying matrix $M$ by $(v,0,0)^{\top}$. $_\Box$
\end{pf}

The next step consists of making explicit the coordinates of the derivative of the unitary vector $\vect{n}$ within frame $(\mathbf{i}_{n},\mathbf{j}_{n},\mathbf{k}_{n})$. This requires two lemmas.

\vspace{0.2cm}

\begin{lemma}
It holds $\frac{d}{dt}\vect{n} = -\frac{1}{R} \left( I - \vect{n} \ \vect{n}\;^{\top} \right) \mathbf{v}$. Moreover, under (\textnormal{\textbf{FOA}}), one has
\begingroup
\footnotesize
\[ \left( \frac{d}{dt}\vect{n} \right)_{(\mathbf{i}_{n},\mathbf{j}_{n},\mathbf{k}_{n})} = \left(\footnotesize \begin{array}{c} 0 \\
\displaystyle -\frac{v}{R} \sin(\gamma - \lambda_1) \\
\displaystyle -\frac{v}{R} \sin(\chi - \lambda_2) \cos(\gamma) \\ \end{array} \right) \]
\endgroup 
\end{lemma}

\begin{pf} Differentiating the definition of $\vect{n}$ one has
\begingroup
\footnotesize
\begin{eqnarray*}
\frac{d}{dt}\vect{n} &= \displaystyle \frac{d}{dt}\left( \frac{\bxi_f - \bxi}{\| \bxi_f - \bxi \|} \right) = -\frac{\mathbf{v}}{R} + \frac{1}{R^2} \frac{(\bxi_f - \bxi) \cdot \mathbf{v}}{R} (\bxi_f - \bxi) \medskip \\
&= \displaystyle  -\frac{\mathbf{v}}{R} + \frac{(\vect{n} \cdot \mathbf{v})}{R} \vect{n} = -\frac{1}{R} \left( I - \vect{n} \ \vect{n}\;^{\top} \right) \mathbf{v}
\end{eqnarray*}
\endgroup
The second expression follows by substituting $(\mathbf{v})_{(\mathbf{i}_{n},\mathbf{j}_{n},\mathbf{k}_{n})}$ within the previous equation.~$_\Box$
\end{pf}

\vspace{0.2cm}

\begin{lemma}
Under (\textnormal{\textbf{FOA}}), we have
\begingroup
$$
\left( \frac{d}{dt}\vect{n} \right)_{(\mathbf{i}_{n},\mathbf{j}_{n},\mathbf{k}_{n})} = \left[ \frac{d}{dt}(M_{(\mathbf{i}_{n},\mathbf{j}_{n},\mathbf{k}_{n})}) \cdot M_{(\mathbf{i}_{n},\mathbf{j}_{n},\mathbf{k}_{n})}^{\top} \right]^{\top} \left(\footnotesize \begin{array}{c} 1 \\
0 \\
0 \\ \end{array} \right)
$$
\endgroup
\end{lemma}

\begin{pf} One has
\begin{eqnarray*}
&\left[ \frac{d}{dt}(M_{(\mathbf{i}_{n},\mathbf{j}_{n},\mathbf{k}_{n})}) \right]^{\top} \left(\scriptsize \begin{array}{c} 1 \\
0 \\
0 \\ \end{array} \right) = \frac{d}{dt}(\vect{n})_{(\mathbf{e}_L, \mathbf{e}_l, \mathbf{e}_r)} \medskip \\
&= \left( \frac{d}{dt} \vect{n} \right)_{(\mathbf{e}_L, \mathbf{e}_l, \mathbf{e}_r)} = M_{(\mathbf{i}_{n},\mathbf{j}_{n},\mathbf{k}_{n})}^{\top} \; \left( \frac{d}{dt}\vect{n} \right)_{(\mathbf{i}_{n},\mathbf{j}_{n},\mathbf{k}_{n}) \ \Box}
\end{eqnarray*}
\end{pf}

\vspace{0.2cm}

Some calculations lead to
$$
\left[ \frac{d}{dt}(M_{(\mathbf{i}_{n},\mathbf{j}_{n},\mathbf{k}_{n})}) \cdot M_{(\mathbf{i}_{n},\mathbf{j}_{n},\mathbf{k}_{n})}^{\top} \right]^{\top} \left(\footnotesize \begin{array}{c} 1 \\
0 \\
0 \\ \end{array} \right) = \left(\footnotesize \begin{array}{c} 0 \\
\dot{\lambda_1} \\
\dot{\lambda_2} \cos(\lambda_1) \\ \end{array} \right)
$$

Gathering together all the previous results, it is straightforward to show the following final result
\begin{equation} \label{lambda1}
\dot{\lambda_1} = \displaystyle -\frac{v}{R} \sin(\gamma - \lambda_1) \; , \quad \dot{\lambda_2} = \displaystyle -\frac{v}{R} \sin(\chi - \lambda_2)
\end{equation}

Relations (\ref{lambda1}) turn out to be useful to obtain initial condition for analytical controls as shown in the followings.

%
%

\subsection{New Approximated Analytical Guidance Law}

\vspace{-0.5cm}

The previous results can be exploited to obtain an approximate solution of (\ref{u1u2App}) and consequently new approximate analytical optimal controls $u_1$, $u_2$ for (\ref{dynamicsSimplyAbs}). We impose $c_m$, $d$ and $\cos(\gamma)$ to be constant along the trajectory (thanks to considerations of Section 3.1). We start considering control $u_1$. From (\ref{u1u2App}), we have
\begin{equation} \label{contU1}
u_1(s) = A e^{\sqrt{\frac{c_m d}{2 \eta}} (s_f - s)} + B  e^{-\sqrt{\frac{c_m d}{2 \eta}} (s_f - s)} - \frac{\cos(\gamma)}{c_m h_r}
\end{equation}

We set $b := \sqrt{\frac{c_m d}{2 \eta}}$. Assumption (\textnormal{\textbf{FOA}}) leads to $R~\cong~(s_f~-~s)$. Plugging (\ref{contU1}) into equation $\gamma' = c_m u_1$ and integrating, one obtains
\begin{equation} \label{intermedium1}
\gamma_f - \gamma(R) = \frac{A c_m}{b}(e^{b R} - 1) - \frac{B c_m}{b}(e^{-b R} - 1)  - \frac{\cos(\gamma)}{h_r}R
\end{equation}

It is straightforward to see that, thanks to (\textnormal{\textbf{FOA}}) and Proposition \ref{PropNeeded}, one has $\dot{R}~=~-v\cos(\gamma - \lambda_1)$. Now, (\ref{lambda1}) is used to differentiate the quantity $R \sin(\gamma - \lambda_1)$. Indeed
\begingroup
\footnotesize
\begin{eqnarray*}
&\frac{d}{dt} (R \sin(\gamma - \lambda_1)) = \dot{R} \sin(\gamma - \lambda_1) + R \cos(\gamma - \lambda_1) (\dot{\gamma} - \dot{\lambda_1}) \medskip \\
&= \dot{R} \sin(\gamma - \lambda_1) + R \cos(\gamma - \lambda_1) \dot{\gamma} + v \sin(\gamma - \lambda_1) \cos(\gamma - \lambda_1) \medskip \\
&= R \cos(\gamma - \lambda_1) v c_m u_1 = -R \dot{R} c_m u_1
\end{eqnarray*}
\endgroup

Then, since under usual assumption $\dot{R} \cong -v$, it follows $(R \sin(\gamma - \lambda_1))' = R c_m u_1$. Integrating this last equation, using (\ref{intermedium1}) and noticing that $R(s_f) = 0$, one has
\begingroup
\footnotesize
\begin{equation*}
\begin{array}{c} \displaystyle -\frac{k_1(R) + k_2(R)}{R^2} \left( R \sin(\gamma - \lambda_1) - \frac{cos(\gamma)}{h_r} \frac{R^2}{2} \right) \medskip \\
\displaystyle -\frac{k_1(R)}{R} \left( \gamma_f - \gamma(R) + \frac{\cos(\gamma)}{h_r} R \right) = c_m u_1(R) + \frac{\cos(\gamma)}{h_r}
\end{array}
\end{equation*}
\endgroup
where $k_1(R)$ and $k_2(R)$ are gain parameters defined as
\begingroup
\footnotesize
\begin{equation*}
\begin{array}{c} \displaystyle k_1(R) = b R \frac{e^{b R} - e^{-b R} - 2 b R}{4 + e^{b R} (b R - 2) - e^{-b R} (b R + 2)} \medskip \\
\displaystyle k_2(R) = b R \frac{e^{b R} (b R - 1) + e^{-b R} (b R + 1)}{4 + e^{b R} (b R - 2) - e^{-b R} (b R + 2)}
\end{array}
\end{equation*}
\endgroup

Since $\sin(\gamma - \lambda_1) \cong \gamma - \lambda_1$, the following guidance law for control $u_1$ is deduced
\begingroup
\begin{equation} \label{lawU1}
\begin{array}{cc} \displaystyle u_1(R) = & -k_1(R) \frac{\gamma_f - \lambda_1(R)}{R c_m} - k_2(R) \frac{\sin(\gamma(R) - \lambda_1(R))}{R c_m} \\ & - k_3(R) \frac{\cos(\gamma(R))}{2 h_r c_m}
\end{array}
\end{equation}
\endgroup
where $k_3(R) = 2 + k_1(R) - k_2(R)$.

With the same argumentation, exploiting again (\ref{lambda1}), the following guidance law for control $u_2$ can be easily derived
\begingroup
\scriptsize
\begin{equation} \label{lawU2}
\begin{array}{c} \displaystyle u_2(R) = -\cos(\gamma(R)) \left(k_1(R) \frac{\chi_f - \lambda_2(R)}{R c_m} + k_2(R) \frac{\sin(\chi(R) - \lambda_2(R))}{R c_m} \right)
\end{array}
\end{equation}
\endgroup

Relations (\ref{lawU1}), (\ref{lawU2}) provide an analytical guess for $(s_f,p_r(0),p_L(0),p_l(0),p_{\gamma}(0),p_{\chi}(0))$ to initialize the shooting method applied to problem (\ref{dynamicsSimplyAbs}). Indeed, they may be chosen as $s_f = R(0)$, $p_{\gamma}(0) = 2 \eta u_1(0)$, $p_{\chi}(0) = 2 \eta \cos(\gamma_0) u_2(0)$ and the guess values of $p_r(0)$, $p_L(0)$, $p_l(0)$ are obtained from (\ref{hamSimply}) and (\ref{dualSys}).

It is interesting to note that (\ref{lawU1}), (\ref{lawU2}) generalize the guidance law proposed in \citep{Lin}.

\section{Solution of the Original System by a Continuation Method}

Since the procedure developed in Section 3 gives us a good guess to initialize the shooting method on a simplified version of (\textbf{OCP}), we proceed now to the analysis of the continuation method that allows to recover an optimal solution of (\textbf{OCP}) starting an homotopy procedure on the simplified problem (\ref{dynamicsSimplyAbs}).

Clearly, several types of different continuations can be implemented. However, we verified that passing from (\ref{dynamicsSimplyAbs}) to (\textbf{OCP}) through a continuation parameter $\lambda_1$ and, after that, changing the final point of the initial scenario directly on (\textbf{OCP}) with another parameter $\lambda_2$ results to be enough to obtain a fast convergence to the optimal solution.

\vspace{-0.3cm}

\subsection{Solution of the Original Problem}

Without loss of generality, we focus on a particular instance of (\textbf{OCP}), whose cost (see (\ref{cost})) takes the form
\begingroup
\footnotesize
$$
\hspace{-4cm} \displaystyle \int_0^{t_f} f^0(t,\vect{x}(t),\vect{u}(t)) \; dt =
$$
$$
\hspace{1.5cm} -\int_{0}^{t_f} \displaystyle \Big[ \frac{f_T}{m} \cos(\alpha) - (d + \eta c_m u^2) v^2 - g \sin(\gamma) \Big] \; dt
$$
\endgroup

The system, which depends on $\lambda_1$, is introduced exploiting the change of variable $w := \ln (v)$ that leaves the formulation consistent with the simplified model (\ref{dynamicsSimplyAbs}) as it is shown hereafter. If we denote $\tilde{\vect{x}}_f$ the final point of the initial test scenario (see Section 3.1) and $\vect{x}_f$ the final point of the desired scenario, we have
\begingroup
\scriptsize
\begin{eqnarray} \label{dynamicsOriginalContinuation}
\min&& \displaystyle -\int_{0}^{t_f} \displaystyle \Big[ \lambda_1 \Big( \frac{f_T}{m e^w} \cos(\alpha) - \frac{g}{e^w} \sin(\gamma) \Big) - (d + \eta c_m u^2) e^w \Big] \; dt \medskip \nonumber \\
\dot{r} =&& e^w \sin(\gamma) \; , \quad \dot{L} = \displaystyle \frac{e^w}{r} \cos(\gamma) \cos(\chi) \; , \quad \dot{l} = \displaystyle \frac{e^w}{r} \frac{\cos(\gamma) \sin(\chi)}{\cos(L)} \medskip \nonumber \\
\dot{w} =&& \displaystyle \lambda_1 \Big( \frac{f_T}{m e^w} \cos(\alpha) - \frac{g}{e^w} \sin(\gamma) \Big) - (d + \eta c_m u^2) e^w \medskip \\
\dot{\gamma} =&& \displaystyle e^w c_m u \cos(\beta) + \lambda_1 \Big( \frac{e^w}{r} \cos(\gamma) + \frac{f_T}{m e^w} \sin(\alpha) \cos(\beta) - \frac{g}{e^w} \cos(\gamma) \Big) \medskip \nonumber \\
\dot{\chi} =&& \displaystyle \frac{e^w c_m}{\cos(\gamma)} u \sin(\beta) + \lambda_1 \Big( \frac{e^w}{r} \cos(\gamma) \sin(\chi) \tan(L) + \frac{f_T}{m e^w} \frac{\sin(\alpha)}{\cos(\gamma)} \sin(\beta) \Big) \bigskip \nonumber
\end{eqnarray}
\endgroup
where the final point takes the form $\vect{x}(t_f) = \tilde{\vect{x}}_f + \lambda_2 (\vect{x}_f - \tilde{\vect{x}}_f)$, $(\lambda_1,\lambda_2) \in [0,1]^2$ and the evolution of the mass is given by $m(t) = m_0 - \lambda_1 \int^t_0 q(s) \; ds$. The continuation algorithm implemented makes $(\lambda_1,\lambda_2)$ converge from $(0,0)$ to $(1,1)$.

As usual, the Hamiltonian of (\ref{dynamicsOriginalContinuation}) takes the form $H = \left \langle \mathbf{p} , \mathbf{f} \right \rangle + f^0$ where we assume that $p^0 \neq 0$ (which is coherent with the simplified formulation (\ref{dynamicsSimplyAbs})). A priori, the Hamiltonian is not uniformly equal to zero because $f_T = f_T(t)$, $q = q(t)$ depend explicitly on $t$. From transversality conditions we obtain
\begingroup
\footnotesize
\begin{equation} \label{transvOriginal}
p_w(t_f) = 0
\end{equation}
\endgroup

Finally, since no constraints are considered on controls $\beta$, $u$, relations $\frac{\partial H}{\partial \beta}~=~0 \; , \ \frac{\partial H}{\partial u} = 0$ hold. From the first one
\begingroup
\footnotesize
\begin{equation} \label{betabeta}
\begin{array}{c} 
\displaystyle 0 = p_{\chi} \frac{c_m e^w u \cos(\beta)}{\cos(\gamma)} - p_{\gamma} c_m e^w u \sin(\beta) \medskip \\
\displaystyle+ \lambda_1 \left( p_{\chi} \frac{f_T}{m e^w} \frac{\sin(\alpha) \cos(\beta)}{\cos(\gamma)} - p_{\gamma} \frac{f_T}{m e^w} \sin(\alpha) \sin(\beta) \right)
\end{array}
\end{equation}
\endgroup

An analytical formulation of $\beta$ is immediately obtained if a first-order approximation $\cos(\alpha) \cong 1$, $\sin(\alpha) \cong \alpha$ is considered (this remains proper if small values of $\alpha_{\max}$ are considered). Indeed, from (\ref{betabeta}) it follows $p_{\chi}\frac{\cos(\beta)}{\cos(\gamma)}~-~p_{\gamma}\sin(\beta)=0$. The same reasoning leads to an explicit relation for $u$.

Variable $w$ is introduced to obtain a guess of $p_w(0)$ for the simplified problem that arises from condition $(\lambda_1,\lambda_2)~=~(0,0)$. In this case, the Hamiltonian is uniformly equal to zero since it does not depend on time. Then, from the dual equations it follows that $\dot{p}_w = -H = 0$ and from (\ref{transvOriginal}) one obtains $p_w(t)~=~0$ in $[0,t_f]$.Then, the procedure achieved in the previous section provides an analytical guess to initialize the sequence of shooting methods within the continuation algorithm that solves the original optimal control problem.

The continuation algorithm is summed up as follows.
\begin{enumerate}
\item[1.] Set $(\lambda_1,\lambda_2) = (0,0)$, $\tilde{\vect{x}}_f = (\tilde{r}_f,\tilde{L}_f,\tilde{l}_f,\tilde{\gamma}_f,\tilde{\chi}_f)$ the final point reached by dynamics (\ref{dynamicsOriginalContinuation}) (equivalent to (\ref{dynamicsSimplyAbs}) at this stage) with control $u = 0$ at $t_f$ and compute the analytical guess for problem (\ref{dynamicsOriginalContinuation}) with initial and final point $\vect{x}_0$, $\tilde{\vect{x}}_f$ respectively;
\item[2.] With $\lambda_2 = 0$ fixed, start a continuation procedure with $\lambda_1~\in~[0,1]$ on the adjoint formulation of (\ref{dynamicsOriginalContinuation}), with initial and final point $\vect{x}_0$, $\tilde{\vect{x}}_f$ respectively, until $\lambda_1~=~1$;
\item[3.] With $\lambda_1 = 1$ fixed, start a continuation procedure with $\lambda_2~\in~[0,1]$ on the adjoint formulation of (\ref{dynamicsOriginalContinuation}), with initial and final point $\vect{x}_0$, $\tilde{\vect{x}}_f+\lambda_2(\vect{x}_f-\tilde{\vect{x}}_f)$ respectively, until $\lambda_2 = 1$.
\end{enumerate}

At point 1 of the previous algorithm $t_f$ must be guessed: it should be chosen such that the assumptions given in Section 3.1 for the integration of (\ref{u1u2App}) are satisfied.

\subsection{Test Cases and Numerical Results}

The previous approach is tested on three test cases. The final scenarios are adopted trying to consider real mission conditions. The following numerical values are chosen: the maximal curvature at altitude zero is chosen as $c_m(0)~=~0.00075 \; \textnormal{m}^{-1}$, the drag coefficient at altitude zero is $d(0)~=~0.00005 \; \textnormal{m}^{-1}$, the aerodynamic efficiency factor is $\eta = 0.442$ and the reference altitude is $h_r = 7500 \; \textnormal{m}$. Moreover, the following functions are chosen to represent respectively the mass flow and the propulsion
$$
q(t) = \left\{
\begin{array}{ll}
      q_0 & t \leq t_{sw} \\
      0 & t > t_{sw}
\end{array} 
\right. \; , \qquad f_T(t) = v_e \cdot q(t)
$$
where $q_0 = 10 \; \textnormal{kg} \cdot \textnormal{s}^{-1}$, $t_{sw} = 20 \; \textnormal{s}$ and $v_e = 1500 \; \textnormal{m} \cdot \textnormal{s}^{-1}$ is the constant fuel injection velocity. Finally, the maximal angle of attack is $\alpha_{\max}~=~\pi / 6$. These conditions are representative of a realistic interceptor scenario. The numerical calculations are computed on a machine Intel(R) Xeon(R) CPU E5-1607 v2 @ 3.00GHz, with 8.00 Gb of RAM.

We impose an initial velocity of $v_0 = 1000 \; \textnormal{m/s}$. If we represent the state by $\vect{y} = (r,L,l,\gamma,\chi)$ (the initial velocity is fixed), the three scenarios are the following (in standard units, where $r_T$ is the Earth radius)
\begin{itemize}
\item[$\bm{S_1}:$] $\vect{y}_0 = (r_T+3000,5454661/r_t,46086/r_t,-\pi/6,0)$ ,

\hspace{-0.45cm} $\vect{y}_f = (r_T+12000,5475000/r_t,42000/r_t,0,\pi/8)$ ;
\item[$\bm{S_2}:$] $\vect{y}_0 = (r_T+3000,5454661/r_t,46086/r_t,\pi/4,0)$ ,

\hspace{-0.45cm} $\vect{y}_f = (r_T+12000,5485000/r_t,36178/r_t,-\pi/4,-\pi/2)$ ;
\item[$\bm{S_3}:$] $\vect{y}_0 = (r_T+3000,5454661/r_t,46086/r_t,0,0)$ ,

\hspace{-0.45cm} $\vect{y}_f = (r_T+3000,5485000/r_t,46086/r_t,0,0)$ .
\end{itemize}

The associated optimal values recovered by the proposed approach are (in standard units) $(v(t_f),t_f)_{\bm{S_1}} = (986.7,24.5)$, $(v(t_f),t_f)_{\bm{S_2}} = (851.6,36.6)$ and $(v(t_f),t_f)_{\bm{S_3}}$ $= (688.8,31.5)$. Using around 400 time steps within the explicit second-order Runge-Kutta solver for the integration of the ODEs inside the shooting function, the simulations take on average $0.27 \; \textnormal{s}$ for $\bm{S_1}$, $0.73 \; \textnormal{s}$ for $\bm{S_2}$ and $2.1 \; \textnormal{s}$ for $\bm{S_3}$. These computing times remain in the order of $1 \; \textnormal{Hz}$ giving the possibility to consider the whole procedure as a good approach for a real-time processing.

The same scenarios were treated also with a non-initialized direct method (AMPL/IPOPT with 200 time steps) \citep{AMPL}. Modifying the boundaries and the initial guess of IPOPT, the three scenarios are solved by the direct method (obtaining the same optimal solutions provided previously) with computational times that vary respectively between $1 \; \textnormal{s}$ and $2 \; \textnormal{s}$ for $\bm{S_1}$, $2 \; \textnormal{s}$ and $5 \; \textnormal{s}$ for $\bm{S_2}$, $5 \; \textnormal{s}$ and $10 \; \textnormal{s}$ for $\bm{S_3}$. Our indirect approach reveals itself to be at least as efficient as a direct method, and sometimes more successful.

It is interesting to note that scenarios $\bm{S_1}$, $\bm{S_2}$ do not require step 3 of the previous algorithm to be solved: imposing $\tilde{\vect{x}}_f~:=~\vect{x}_f$ within step 1 , the initial guess is good enough to find the solution of (\textbf{OCP}) with a continuation on $\lambda_1$ only. From this, 3 iterations of the continuation method on $\lambda_1$ are required to solve $\bm{S_1}$, 10 iterations on $\lambda_1$ for $\bm{S_2}$ while 18 iterations on $\lambda_1$ and 6 iterations on $\lambda_2$ are needed to solve $\bm{S_3}$.

The trajectories and the optimal controls $u_1~:=~u\cos(\beta)$, $u_2~:=~u\sin(\beta)$ are shown in Figure 1. It is interesting to note that, even if no explicit constraint are considered on both the controls, they turn out to be bounded in $[-1,1]$. They are forced to be small due to the presence of $u^2$ inside the cost function.

\begin{figure*}[htpb]
\centering
\includegraphics[width=1\textwidth]{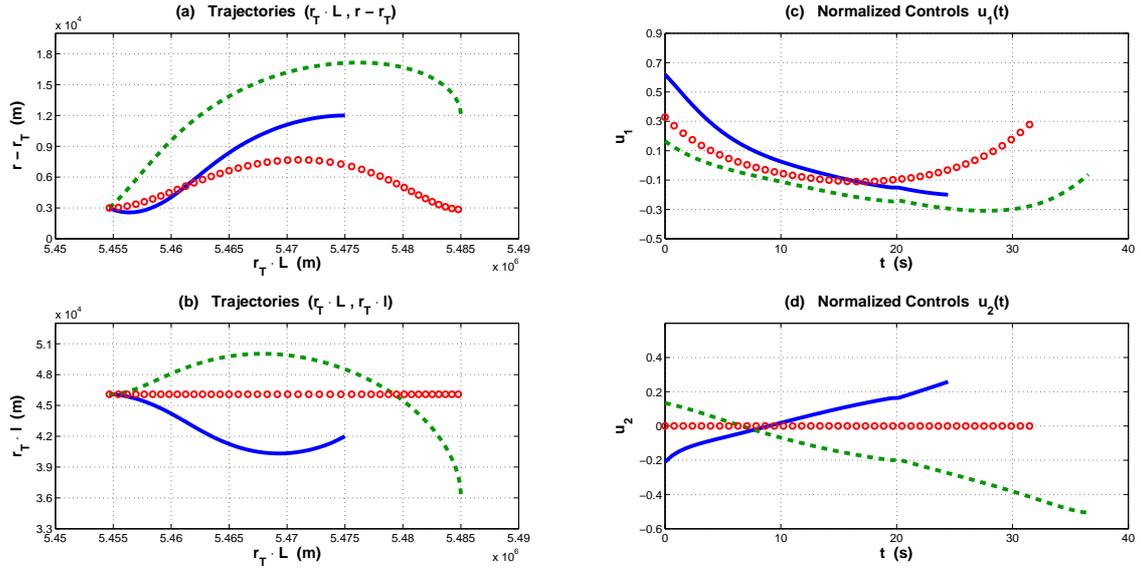}
\caption{Optimal trajectories and controls $u_1$, $u_2$. The solid line represents $\bm{S_1}$, while the dashed and the empty-dot lines represent respectively $\bm{S_2}$ and $\bm{S_3}$.}
\label{test2}
\end{figure*}

\section{Conclusions and Perspectives}

In this paper we solve an optimal control problem for endo-atmospheric launch vehicle systems using indirect methods. The problem is solved mixing an analytical initialization with a continuation method on the dynamics. First, a simplified version of the original problem is solved. Then, a continuation method which recovers the complete optimal solution starts: every shooting method within the iterative procedure is initialized by the approximate solution found at the previous stage. This approach provides numerical precision and high computational. Moreover, numerical simulations show that this approach may compete with usual direct methods.

Clearly, some improvements can be brought. A first remark concerns the model. Indeed, the velocity of the wind and stabilizing components within the cost such as $\sigma u^2$ have not been considered in the previous analysis. However, such terms can be added to the formulation exploiting a new stage of continuation on them. Secondly, in this contribution, neither control constraints nor state constraints were taken into account. It is not hard to consider control constraints (a detailed geometric study \citep{mio} shows that neither abnormal nor singular arcs exist for this problem) while, as it is known, state constraints result to be more challenging (one could, for example, penalize these constraints within the cost function \citep{NewEmmanuel}). Finally, as in previous sections, dynamics (\ref{dynamics}) reveals some singularities due to eulerian coordinates. In order to deal with this flaw, an inline change of coordinates can be considered during the numerical procedure. This strategy does not increase the computational effort and it allows to solve many more scenarios \citep{mio}.

All these precious improvements will be reported in \citep{mio}.

\bibliography{ifacconf}             

\end{document}